\magnification=1200
\hsize =15true cm
\baselineskip = 16pt
\vsize= 22true cm

\centerline { \bf A Relation between $\Gamma$-Convergence of Functionals }

\centerline {\bf and their Associated Gradient Flows* }

\vskip 0.3cm

\centerline { \bf JIAN Huaiyu }

\centerline {  Department of Applied Mathematics, Tsinghua University,
Beijing 100084 }

\vskip 0.3cm

{\bf Abstrract } \ \  E. De Giorgi conjectured in 1979 that if a sequence
of functionals converges in the sense of $\Gamma$-convergence to a limiting
functionals, then the corresponding gradient flows will converge as well
after changing timescale appropriately. In this paper, we will show this
conjecture holds true for a rather wide kind of functionals.


\vskip 0.3cm

 In 1979, E. De Giorgi [1] asked if there was a general relation between
 $\Gamma$-convergence of functionals and convergence of solutions to the
 associated parabolic equations. He further conjectured in the same paper
 that when a sequence of functionals converges in the sense of $\Gamma$-convergence
  to a limiting functional, then the corresponding gradient flows will
  converge as well (maybe after an approriate change of timescale). Also see
  [2; P.216] and [3; P.507]. Although there is no any result, up to the author's
  knowledge, confirming this conjecture, it was supported by the results of
  Bronsard and Kohn in [2], and Owen, Rubinstein and Sternberg in [3], respectively,
  where they studied the singular limit of Ginzburg-Landau dynamics (up to a
  $\varepsilon$-scaling time):
  $$u_t - \varepsilon ^2 \Delta u +u^3 - u =0  \eqno (0.1) $$
  which are  the gradient flows of the following  functionals:
 $$f_{\varepsilon}(u) = \int_{\Omega} ({1\over 4} (u^2 -1)^2 +{\varepsilon^2\over
 2}|\nabla u|^2 )dx, \ \ \Omega \subset R^n . \eqno (0.2) $$
 The $\Gamma$-limit of functionals (0.2) as $\varepsilon \to 0^+$ was derived by Modica in [4]. Combination
 of the results in [2, 3] with ones in [4] suggests that De Giorgi's conjecture
 should be answered positively, at least for some special functionals.

\noindent $\overline{* Supported \ \ by \ \ N N S F \ \ of \ \ China  \ \ (Grant
\ \ No. 19701018) \ \ \ \ \ \ }$

 In this paper, we will confirm De Giorgi's  conjecture for a rather wide kind of
 functionals. Precisely speaking, we will establish a relation between
 $\Gamma$-convergence of functionals and the convergence of their parabolic-minima.
 The $\Gamma$-convergence may be proved by similar arguements in [5, 6].
 Furthermore, we discover that parabolic-minima of rather many functionals are
 nothing but the solutions to the gradient flows of the corresponding functionals.

\

\

\noindent {\bf 1 \ \ Main results }

\vskip 0.3cm

We begin with the following assumptions and notations:
 $\Omega$ denotes a bounded open set  in $ R^n$ , $p>1$, $T>0$ and $m$ is
 a positive integer. Let
$\Omega _T=\Omega\times (0,T)$,
$$V_p(\Omega _T,m)=L^P([0 , T],W^{1,p}(\Omega , R^m)), \ \
V_p^0(\Omega _T ,m)=L^p([0,T],W_0^{1,p}(\Omega , R^m))  ,\eqno (1.1) $$ and
$$ Du(x,t)=\nabla u (x, t)=\bigl ({\partial u^i(x,t) \over \partial x_j } \Bigr ) \quad
(1 \leq i \leq m , 1 \leq j \leq n )  \eqno (1.2) $$ for a vector valued function $u$.

Suppose $\Phi \colon R^m \to R^m,$
 and
$f(x,t,u, \lambda ) \colon \Omega _T \times R^m \times R^{mn} \to R $ such that
for each $v \in V_p(\Omega _T ,m)$
$$\Phi (v) \in L^1(\Omega_T ,R^m)  \ \ and \ \   f(x,t,v,Dv) \in
L^1(\Omega _T). \eqno (1.3) $$
\noindent Consider the parabolic functional
$$ F(v, \Omega _T)=\int _{\Omega_T}f(x,t,v,Dv)dxdt \quad v \in V_p(\Omega_T, m).
\eqno (1.4)  $$
Following the idea of the papers [7], we introduce the definition
of parabolic-minima.

  \vskip 0.3cm

{\bf Definition 1.1.}  \ \   {\sl Assume that $\Phi $ and $f$ satisfy (1.3),
 $u_0(x) \in L^1(\Omega ,R^m). $   A function
$u \in V_p(\Omega _T, m)$ is called a parabolic-minimum of $F$(defined by (1.4))
with respect to
the function couple $(\Phi , u_0 )$
if for all $\eta  \in C^\infty
([0,T], C_0^\infty
(\Omega , R^m)) $ with $ \eta (\cdot,T)=0,$
$$ - \int _{\Omega _T}\Phi (u){\partial \eta \over \partial t }dxdt +
F(u,\Omega _T)
\le F(u- \eta , \Omega _T) + \int _\Omega u_0\eta (x,0)dx . \eqno (1.5) $$}

{\bf Definition 1.2.} \ \ {\sl
 $\tau $ is called as the {\bf sw}-topology of $V_p(\Omega _T,m),$
 if $ v^{\varepsilon}$ converges to $v$ in $V_p (\Omega_T , m)$ with the
 topology $\tau$ if
and only if $ v^\varepsilon \to v$ strongly in $L^p(\Omega_T , R^m)$  and
 $ Dv^\varepsilon \rightharpoonup
Dv$  weakly in  $L^p(\Omega _T , R^{mn}).$
 We denote this convergence by
 $v^\varepsilon \buildrel \tau \over\rightarrow v.$}

\vskip 0.3cm

Now consider a sequence of functionals defined in $V_p (\Omega _T ,m)$
by
$$F^\varepsilon (v, \Omega _T) =\int _{\Omega _T}f^\varepsilon (x,t,v,Dv)dxdt ,
(\varepsilon \rightarrow 0 ),  \eqno (1.6) $$
where each $f^{\varepsilon} \colon \Omega_T \times R^m \times R^{mn} \rightarrow R $ is a
{\sl Caratheodory } function satsfying
$$ 0 \le f^\varepsilon (x,t,u,\lambda ) \le C (1+ |u|^p+|\lambda |^p )
 \eqno (1.7)$$ and
$$\eqalignno {|f^\varepsilon (x,t,u_1,\lambda _1)-f^\varepsilon (x,t,u_2
,\lambda _2)|
\le & C (|u_1-u_2|^\alpha +|\lambda _1 - \lambda _2 |^\alpha )(1+|u_1|^{p- \alpha }  \cr
&  +|u_2|^{p- \alpha }+| \lambda _1|^{p- \alpha }+| \lambda _2|^{p- \alpha })
&(1.8)   \cr} $$
\noindent for some constants $ C > 0$    and   $ \alpha \in (0, 1).$

\vskip 0.3cm

The main result of this paper is the following theorem where we refer to [6, 8]
for the definition of $\Gamma$-convergence.

{\bf Theorem 1.1.} \ \  {\sl Suppose that the hypotheses (1.6), (1.7) and  (1.8)
hold true  and that $F^\varepsilon$  $\Gamma $-converges to $F$ with {\bf sw}-topology,
i.e.,
$$ \Gamma (\tau ) \lim_{\varepsilon  \to 0 }F^\varepsilon (v,
\Omega _T) = F(v, \Omega _T) , \forall v \in V_p(\Omega _T, m),\eqno (1.9) $$
\noindent where $\tau$ is the {\bf sw}-topology of $ V_p(\Omega _T, m).$
 If for each $\varepsilon > 0, $  $ u^\varepsilon
\in V_p^0(\Omega _T, m) $ is a parabolic-minimum of $F^\varepsilon $ with respect
to $ (\Phi  , u_0^\varepsilon ) $ such that as $\varepsilon \to 0 ,$
$$u_0^\varepsilon  \rightharpoonup u_0 \ \  and \ \ \Phi
(u_0^\varepsilon)
 \rightharpoonup \Phi (u_0)  \ \ weakly \ \ in \ \ L^1(\Omega _T , R^m)  \eqno
 (1.10)$$ and
$$ u^\varepsilon \buildrel \tau \over
\rightarrow u \in V_p^0(\Omega _T) , \ \ and \ \ \partial_t \Phi (u^\varepsilon )
\rightharpoonup \partial_t \Phi (u)  \ \ weakly \ \ in \ \ L^q(\Omega _T,R^m)
\eqno (1.11) $$
with $q={p \over p-1},$
 then $u$ is a parabolic-minimum of $F(u)$ with
respect to $ ( \Phi ,u_0 ) $. }

\vskip 0.4cm

We will prove this theorem in section 3, while in next section, we will study the
equivalence of parabolic-minima with some parabolic systems and discuss the
 justification for assumptions (1.10) and (1.11). We would like to point out that
 assumption (1.9) may checked by arguements similar to those in [5, 6, 8].
 In the case of $f^{\varepsilon}(x, t, u, \cdot )=f({x \over \varepsilon },{t
 \over \varepsilon },{u\over \varepsilon},\cdot)$ with $f$ being periodic in
 the first three variables, it was proved in [9].

\

\

\noindent {\bf 2 \ \  Parabolic-minima and parabolic equations }

\vskip 0.3cm

{\bf Lemma 2.1.} \ \   {\sl Suppose $ f \in C^1 $ in $ v $ and $ \lambda $ such that
$$ |f(x,t,v,\lambda )| + |f_v(x,t,v, \lambda )|+|f_\lambda (x,t,v, \lambda )|
\le C (1+|v|^p +|\lambda |^p)  .  \eqno (2.1) $$
\noindent If $ u \in V_p^0( \Omega _T)$ is a parabolic minimum of the functional $F(v,\Omega_T)$
given by (1.3)-(1.4) with respect to $(\Phi ,u_0)$, then u is a weak solution
to its gradient flow which is the
 following initial-boundary valued problem:}
$$\cases { \Phi (u) \in L^1(\Omega _T , R^m) , \ \ and  \ \ for
 \ \  i=1,2, \cdots , m  \cr
{\partial \Phi^i (u)  \over \partial t } - {\partial f_{\lambda _\alpha ^i}
\over \partial x_\alpha }(x,t,u,Du)+f_{u^i}(x,t,u,Du) =0 \ \ in \ \  \Omega _T, \cr
u=0 \quad  on \quad \Sigma = \partial \Omega \times (0,T) , \quad  \Phi (u(x,0))=u_0(x)  \cr }
\eqno (2.2) $$

{\bf \sl Proof.}  \ \ For each $ \eta \in C^\infty ([0,T], C_0^\infty (\Omega , R^m))$
with $\eta (\cdot ,T)=0 $ and any $ h \in R\backslash \lbrace 0 \rbrace ,$  replace
$\eta$ by $h \eta $ in (1.5), use the mean value formula and then divide by h.
Letting $ h \longrightarrow 0^+ $ and $ h \longrightarrow 0^- $ respectively,
we  obtain the desired result.

\vskip 0.3cm
{\bf  Remark  2.1.} \ \ If $F(v, \Omega _T )$ has a parabolic minimum and (2.2)
has a unique solution, then lemma 2.1 implies that the solution must be the
 parabolic-minimum. But we don't know the existence for the parabolic-minimum.
 Nevertheless, the following examlples show that parabolic-minima of some
 functionals are nothing but weak solutions to their corresponding
 gradient flow equations.

\vskip 0.3cm

{\bf Example 2.1.} \ \ For each parameter $\varepsilon ,$ let
$  \bigl [a_{ij}^
\varepsilon (x, t) \bigr ]  $
  be a symmetric, measurable and uniformly
bounded and positive-definite matrix function in $\Omega _T .$
Assume that $\Phi =(\phi ^1 , \cdots , \phi ^m)$ is a map from $R^m$ to
itself with monotonic components.  Consider the asymptotic
behaviour of weak solutions to the equations of general Newtonian filtration :
$$  \cases { {\partial \phi ^k (u) \over \partial t} =
{\partial \over \partial x_i}(a_{ij}^
\varepsilon (x,t){\partial u^k \over \partial x_j} )+h^k_\varepsilon (x,t)
\ \  in  \ \  \Omega _T  , k=1,2, \cdots , m   \cr
 u=0 \qquad on \ \  \partial \Omega \times (0,T)   \cr
u(x,0) = u_0(x) .  \cr }  \eqno (2.3) $$
It is well-known that under some additional assumptions, for instance,
$h_\varepsilon =(h^1_\varepsilon , \cdots , h^m_\varepsilon) \in L^2
(\Omega _T ,R^m)$ and all $\phi ^k $ have boumded derivatives, the solutions
$u^\varepsilon $ to (2.3) belong to $V_2(\Omega _T ,m)$ and
$$\| \partial_t \Phi (u^\varepsilon )\|_{L^2(\Omega _T)} + \|u^\varepsilon \|
_{V_2(\Omega _T)}\leq C   \eqno (2.4)$$
for some constant $C$ independent of $\varepsilon .$  See [10, 11]. Using
this result, lemma B below and the compact imbedding theorem, we have
a sequence of
$u^\varepsilon $ such that it satisfies (1.10) and (1.11) for $p=2 .$ Moreover, we have
 the following conclusion.

\vskip 0.3cm

{\bf Lemma 2.2.} \ \ For each $\varepsilon ,$ the solution $u^\varepsilon $ to
(2.3 ) is a parabolic-minimum of the functional
 $$F_1^\varepsilon (v)= \int _{\Omega _T}({1 \over 2}
 a_{ij}^\varepsilon (x,t)
{\partial v^k  \over \partial x_i }{\partial v^k \over \partial x_j}- h^k_\varepsilon
v^k)dxdt
  \ \ \ v\in V_2(\Omega _T ,m) \eqno (2.5)$$
with respect to $(\Phi , u_0 ) .$

{\bf \sl Proof.}  \ \ For simplicity, we denote $u^\varepsilon$ by $u.$ Since u is
a weak solution to (2.3), we see that
 for each $\eta  \in C^\infty
([0,T], C_0^\infty
(\Omega , R^m)) $ with $ \eta (\cdot,T)=0,$
$$\eqalignno{ - \int _{\Omega _T}\Phi (u){\partial \eta \over \partial t }dxdt
+ F^\varepsilon (u) = &
{1\over 2}
\int_{\Omega _T}
 a_{ij}^\varepsilon
{\partial u  \over \partial x_i }{\partial u
 \over \partial x_j}dx dt -
\int_{\Omega _T}
 a_{ij}^\varepsilon
{\partial u  \over \partial x_i }{\partial \eta
 \over \partial x_j}dx dt  \cr
& - \int _{\Omega _T} h_\varepsilon (u-\eta) dxdt
 + \int _\Omega u_0\eta (x,0)dx . & (2.6) \cr } $$
As
$  \bigl [a_{ij}^
\varepsilon (x, t) \bigr ]  $
is positive-definite, the sum of first two terms on the right hand in (2.6)
is no larger than
$${1\over 2}
\int_{\Omega _T}
 a_{ij}^\varepsilon
{\partial (u-\eta )  \over \partial x_i }{\partial (u-\eta)
 \over \partial x_j}dx dt .$$
Thus, (2.6) turns to
$$ - \int _{\Omega _T}\Phi (u){\partial \eta \over \partial t }dxdt +   F_1^\varepsilon
(u)
\leq F_1^\varepsilon (u-\eta )
 + \int _\Omega u_0\eta (x,0)dx ,  $$
 as desired.

\vskip 0.3cm

 {\bf Example 2.2.}
\ \ Suppose  that $a^\varepsilon (x,t) $ is a family of measurable functions
which are positive and bounded uniformly in parameter $\varepsilon .$
Let $ p>1 $ and $\Phi $ be the same as in example 2.1. If for each
$\varepsilon ,$
$u^\varepsilon $
is the weak solution to the following P-Laplace equation :
$$  \cases  {\partial_t \phi ^k  - \partial_{x_\alpha} (a^\varepsilon (x,t)
|\nabla u |^{p-2} \partial_{x_\alpha}  u )=0  \ \  in  \ \  \Omega _T
 , k=1,2, \cdots ,  m   \cr
 u=0 \ \   on   \ \  \partial \Omega \times (0,T) , \quad u(x,0)=u_0(x) . \cr }
 \eqno (2.7) $$
 then $u^\varepsilon  $
is a parabolic minimum of
$$ F_2^\varepsilon (v) =
{1 \over p } \int _{\Omega _T}a^\varepsilon (x,t) | \nabla
v |^p dxdt   . $$
The proof is similar to the arguements used in proving lemma 2.2. We omit
the details.

\vskip 0.3cm

{\bf Example 2.3.} \ \ Let $u_\varepsilon$ be the solutions to the Cauchy
problem of equation (0.1) with initial data $u^\varepsilon_0 $ satisfying
$u^\varepsilon_0 \geq 1 .$ Obviously, $|u_\varepsilon (x, t) |\geq 1 ,$
for all $(x,t)\in R^n \times (0, \infty ) ,$ by the maximum principle.

Since for each $\eta  \in C^\infty
([0,T], C_0^\infty
(\Omega , R^m)) $ with $ \eta (\cdot,T)=0,$
 we obtain that
$$\eqalignno{ - \int_0^T & \int _{R^n}  u_\varepsilon {\partial \eta \over \partial t }dxdt
 +{1\over 2}
\int_0^T\int_{R^n}\bigl [\varepsilon^2 |\nabla u_\varepsilon | +(u_\varepsilon^2 -1)^2
\bigr ]dxdt
- \int _{R^n} u^\varepsilon_0 \eta (x,0)dx \cr
& ={\varepsilon ^2\over 2} \int_0^T \int_{R^n}\big [ |\nabla (u_\varepsilon -\eta )
|^2 -|\nabla \eta |^2 \bigr ]dxdt \cr
& \ \ \ +{1\over 2}\int_0^T \int_{R^n}\bigl [
(u_\varepsilon ^2 -1)((u_\varepsilon -\eta )^2 -1)-\eta ^2 (u_\varepsilon ^2
-1)\big ] dxdt \cr
& \leq {\varepsilon ^2\over 2} \int_0^T \int_{R^n} |\nabla (u_\varepsilon -\eta )
|^2 dxdt \cr
& \ \ \  +{1\over 4}\int_0^T \int_{R^n}
(u_\varepsilon ^2 -1)^2 dxdt+{1\over 4}\int_0^T\int_{R^n}\big [ (u_\varepsilon
- \eta )^2 -1 \big ]^2 dxdt, \cr}$$ where we have used the fact that
$u_\varepsilon ^2 -1\geq 0$ and Young's inequality.
This immediately implies that $u_\varepsilon$ is a parabolic-minimum of
functionals of type (0.2), i.e.,$u_\varepsilon$ is a paraboilc-minimum of
 $$F_3^\varepsilon (u)= \int_0^T\int_{R^n} \bigl ({1\over 4} (u^2 -1)^2 +{\varepsilon^2\over
 2}|\nabla u|^2 )dx  $$ with respect to $( I, u_0^\varepsilon )$, where $I$ is the identity
 map.

\vskip 0.3cm

 {\bf  Remark 2.2.} \ \  For the initial-boundary value problem of (0.1),
 we also have a similar
 conclusion.

 \

 \

\noindent  { \bf 3 \ \   A proof of theorem 1.1 }

\vskip 0.3cm

To prove theorem 1.1, we need two well-known results.

{\bf Lemma A.}  \ \  {\sl Suppose that $ f \colon R\times R \to \bar R$, then
there exists a function $\delta \colon \varepsilon \to \delta (\varepsilon )$
such that $\varepsilon \to 0 $ implies $ \delta (\varepsilon ) \to 0 $ and
$$\lim_{\varepsilon \to 0 }f (\delta (\varepsilon ),\varepsilon )
=\lim_{\delta \to 0 } \lim_{\varepsilon \to 0 }
f(\delta , \varepsilon )   $$ }

{\bf \sl Proof. } See  [5; P.32-33].

\vskip 0.3cm

{\bf Lemma B. } \ \    {\sl If $p>1 $,  then

(1) $\{u_h \}$ converges to $u$ with respect to {\bf sw}-topology of
$V_p(\Omega_T ,m)$ if and only if

\ \ $u_h\rightarrow u $ strongly in $L^p(\Omega_T ,R^m)$
and $ \sup_h\|\nabla u_h\|_{L^p(\Omega_T, R^{mn})} <  \infty ; $

(2) the fact that $\{u_h\}$ converges to $u$ with respect to  {\bf sw}-topology
of  $V_p(\Omega _T , m)$ implies that $ \sup_h\|u_h\|_{V_p(\Omega_T ,m)}
<  \infty  .$ }

{\bf \sl Proof.} See section 2 of chapter 1 in [12].

\

\

Now we are in the position to prove theorem 1.1.  Let $\tau$ be the
 {\bf sw}-topology. For simiplicity, denote $V_p(\Omega _T ,m)$ by
 $V_p(\Omega _T ).$

 Fix $\eta \in C^\infty ([0,T],C_0^\infty (\Omega , R^m))$ arbitrarily.
According to assumption ( 1.9) and the definition of $\Gamma$-convergence[8],
 we can choose
 $ \lbrace w^\varepsilon \rbrace \subset
V_p(\Omega _T)$ satisfying
$$ w^\varepsilon \buildrel \tau \over \rightarrow u-\eta \eqno (3.1) $$
suct that
$$ F(u- \eta ,\Omega _T) = \lim_{\varepsilon \to 0}F^\varepsilon (w^\varepsilon ,\Omega _T) .
\eqno (3.2) $$
\noindent Since $ u-\eta \in V_p^0(\Omega _T),$  we can assert that there exist
a sequence $\{v^\varepsilon \}  \subset  V_p^0(\Omega_T)$  such that
$$v^\varepsilon \buildrel \tau \over \rightarrow u-\eta  \eqno (3.3) $$
and
$$\lim_{\varepsilon \to 0}F^\varepsilon (v^\varepsilon , \Omega_T)=F(u-\eta ,
\Omega_T) . \eqno (3.4)  $$

Indeed, for each $ Q_0 \subset \subset
\Omega_T ,$  let $R=2^{-1} dist(Q_0, \partial \Omega_T).$  For any $\delta
\in (0, 1) $, define
$$ Q_i =\{ x \in \Omega_T \colon   dist(x, Q_0) <i \delta R\}, i=1,2,\cdots
 [\delta^{-1}] . $$
Choose $\phi_i \in C_0^\infty (\Omega_T)$ such that $0\le \phi_i \le 1 , \phi_i
=1 \ \  in  \ \  Q_{i-1}, \phi_i=0 \ \ in \ \ \Omega_T \backslash Q_i$ and
$|D\phi_i| \le C(\delta , n, R) .$  let
$$ u^{\varepsilon, i} = u-\eta +\phi_i (w^\varepsilon -u+\eta).$$
By (1.7) we have
$$\eqalignno {F^\varepsilon(u^{\varepsilon ,i}, \Omega_T) & =F^\varepsilon(
u^{\varepsilon, i}, Q_{i-1})+F^\varepsilon(u^{\varepsilon, i}, Q_i\backslash
Q_{i-1})+F^\varepsilon(u-\eta , \Omega_T \backslash Q_i) \cr
& \le F^\varepsilon(w^\varepsilon, \Omega_T)+C\bigl [\|u-\eta\|_{V_p(\Omega_T\backslash
Q_0)}^p +\|w^\varepsilon\|_{V_p(Q_i\backslash Q_{i-1})}^p \cr
& +C(\delta, R)\|w^\varepsilon -u+\eta\|_{L^p(\Omega_T)}^p+|\Omega_T\backslash
Q_0| \bigr ] .  & (3.5) \cr } $$
Obviously, we can  find $i(\delta) \in \{1,2, \cdots , [\delta^{-1}] \}$ such that
$$ F^\varepsilon (u^{\varepsilon, i(\delta)}, \Omega_T)= \min \{F^\varepsilon(
u^{\varepsilon, i},\Omega_T) \colon i=1,2, \cdots [\delta^{-1}] \} . $$
Furthermore, (3.1) implies
$$ u^{\varepsilon,i(\delta )} \buildrel \tau \over \rightarrow u-\eta
\ \ for \ \ each \ \ \delta . \eqno (3.6) $$
Summing (3.5) for $i$ from 1 to $[\delta ^{-1} ]$, we arrive at the
estimate
$$\eqalignno{F^\varepsilon(u^{\varepsilon,i(\delta)}, \Omega_T)\le & F^\varepsilon
(w^\varepsilon, \Omega_T)+C\bigl [\|u-\eta\|_{V_p(\Omega_T\backslash Q_0)}^p+{[\delta
^{-1}]}^{-1}\|w^\varepsilon\|_{V_P(\Omega_T)}^p \cr
& + C(\delta ,R) \|w^\varepsilon -u+\eta \|_{L^p(\Omega_T)}^p +|\Omega_T
 \backslash Q_0| \bigr ] . \cr }$$
Applying this estimate, (3.6), (1.9), the definition of $\Gamma$-convergence[8],
 (3.2) and lemma B,
we obtain that
$$\eqalignno{F(u-\eta, \Omega_T) & \le \liminf_{\delta \to 0}\liminf_{\varepsilon
\to 0}F^\varepsilon
(u^{\varepsilon, i(\delta)},\Omega_T) \cr
& \le \limsup_{\delta \to 0}\limsup_{\varepsilon
\to 0}F^\varepsilon
(u^{\varepsilon, i(\delta)},\Omega_T) \cr
& \le F(u-\eta, \Omega_T)+C\bigl ( \|u-\eta\|_{V_p(
\Omega_T\backslash Q_0)}^p +|\Omega_T\backslash Q_0|\bigr ) . \cr }$$
 Letting $Q_0 \to \Omega_T$, we have
$$\lim_{\delta \to 0} \lim_{\varepsilon \to 0}F^\varepsilon(u^{\varepsilon,i(\delta)},
\Omega_T)= F(u-\eta ,\Omega_T) .  \eqno (3.7)$$
Now define
$$\eqalignno{f(\delta,\varepsilon )= & |F^\varepsilon(u^{\varepsilon,i(\delta)}
,\Omega _T)-F(u-\eta , \Omega _T)| \cr
& +\|u^{\varepsilon,i(\delta)}-u+\eta\|_{L^p(\Omega_T)}+\|D\phi_{i(\delta)}(
w^\varepsilon-u+\eta)\|_{L^p(\Omega_T) }. \cr } $$
Then (3.6), (3.7) and lemma B implies that
$$\lim_{\delta \to 0}\lim_{\varepsilon \to 0} f(\delta ,\varepsilon) =0 .$$
By virture of lemma A, we conclude that there exists a function
 $\delta(\varepsilon)$  such that
 $$\lim_{\varepsilon
\to 0}f(\delta(\varepsilon), \varepsilon) =0 .$$
 Therefore,  setting $i(\varepsilon)=i(
\delta(\varepsilon)) $  and $  v^\varepsilon =u^{\varepsilon, i(\varepsilon)}$,
we  see that $\{v^\varepsilon \}\subset V_p^0 (\Omega _T)$ satisfy (3.4).
Moreover, we have that as $\varepsilon \to 0$,
$$\|v^\varepsilon -u+\eta \|_{L^p(\Omega_T)} \to 0  \quad and \quad
\|D\phi_{i(\varepsilon)}
(w^\varepsilon -u+\eta) \|_{L^p(\Omega_T)} \to 0, $$
which, together with lemma B, implies that $\{v^\varepsilon \}$ satisfy
(3.3).

\vskip 0.3cm

Now take $\{u^\varepsilon \}$ as in theorem 1.1. Then for   each $\varepsilon $,
 $u^\varepsilon - v^\varepsilon  \in V_p^0(\Omega _T)$. Thus we can choose
 $ \eta ^\varepsilon \in C^\infty ([0, T], C_0^\infty (\Omega , R^m))$
such that
$$\lim_{\varepsilon \to 0} \| \eta ^\varepsilon - (u^\varepsilon -
 v^\varepsilon ) \| _{V_p(\Omega _T)}=0 .$$
\noindent Furthermore, we may assume
$$ \eta ^\varepsilon (x,0) = \eta (x,0) , \quad \eta ^\varepsilon (x,T) =0 .
\eqno (3.8) $$
By (1.8), we get that
$$|F^\varepsilon (v^\varepsilon , \Omega _T)-F^\varepsilon (u^\varepsilon  -
\eta ^\varepsilon , \Omega _T)|  \le C(m,n,p)  \times $$
$$\| u^\varepsilon - v^\varepsilon  - \eta ^\varepsilon \|_{V_p(\Omega _T)}^ \alpha \Bigl (1+
\|\eta ^\varepsilon \|_{V_p(\Omega _T)} + \|u^\varepsilon \|_{V_p(\Omega _T)}
+\| v^\varepsilon \|_{V_p(\Omega _T)} \Bigr )^{p- \alpha }.$$
\noindent Observing that the sequence $\{u^\varepsilon\} \ \  and \ \  \{v^\varepsilon\} $
 are bounded (see lemma B ), we obtain
$$\lim_{\varepsilon \to 0 }F^\varepsilon (v^\varepsilon , \Omega _T) = \lim_{
\varepsilon \to 0 } F^\varepsilon (u^\varepsilon - \eta ^\varepsilon , \Omega _T)
. \eqno (3.9) $$
\noindent Since for each $\varepsilon$, $u^\varepsilon $ is a parabolic-minimum  of $F^\varepsilon $,
it follows  from (3.4), (3.9) and (1.10) that
$$ F(u- \eta ,\Omega_T) \ge \liminf_{\varepsilon \to 0 }\Bigl ( F^\varepsilon (u^\varepsilon ,
\Omega _T) - \int _{\Omega _T}\Phi  (u^\varepsilon ){\partial \eta
^\varepsilon \over \partial t}dxdt \Bigr )- \int _\Omega u_0(x)\eta (x,0) dx  . \eqno (3.10) $$
\noindent By (1.11) and (3.3),
 $u^\varepsilon - v^\varepsilon \buildrel \tau \over\rightarrow
\eta$,
 so $\eta^\varepsilon \buildrel \tau \over\rightarrow \eta . $
Therefore, again by  assumption (1.11),
 we have  $$ \lim_{\varepsilon \to 0} \int_{\Omega_T} \Phi
(u^\varepsilon) \partial_t \eta^\varepsilon dxdt
=\int_{\Omega_T}\Phi (u)\partial_t \eta dxdt .$$
\noindent Thus (3.10) yields $$ F(u-\eta ,\Omega _T) \ge \liminf_{\varepsilon
\to 0 } F^\varepsilon (u^\varepsilon ,\Omega _T) - \int _{\Omega _T}\Phi (u)
{\partial \eta \over  \partial t }dxdt - \int _\Omega u_0(x) \eta (x,0)dx .$$
\noindent Noting
 $u^\varepsilon  \buildrel \tau \over\rightarrow  u$ (by (1.11)) and using (1.9),
 definitions 1.1 and the definition of $\Gamma$-convergence[8],  we have completed the proof of theorem 1.1.

\

\

\

\

\

\noindent {\bf References }

\item{1}  De Giorgi E.  New problems in $\Gamma$-convergence and G-convergence
in Free Boundary Problems. In: Nazeson I, Materson A, eds.  Proc. of seminar held
in Pavia. Pivia, September-October
1979, Rome, Francesco Severi,  1980, Vol. II 183-194

\item{2} Bronsard L,  Kohn R.  Motion by mean curvature as the singular limit
of Ginzburg-Landau Dynamics. {\sl J Differ Eqn}, 1991, 90:  211-237

\item{3}  Owen N,  Rubinstein J,  Sternberg P. Minimizers and gradient
flows for singularly perturbed bi-stable potentials with a Dirichlet
condition. {\sl Proc R Soc Lond}, 1990,429(A):  505-532

\item{4}  Modica L. The gradient theory of phase transitions and the minimal
interface criterion. {\sl Arch Rat Mech Anal}, 1987, 98:  123-142

\item{5}  Attouch H.  Variational convergence for functions and operators.
Pitman:
 Appl Math Series, 1984. 23-198

\item{6}  Dal Maso G. An introduction to $\Gamma$-convergence. Boston Basel:
Birkh\"auser, 1993. 38-94.

\item {7}  Wieser W.  Parabolic Q-minimum and minimal solutions to
varational flow. {\sl Manu Math}, 1987,    59: 63-107

\item {8}  Jian H. $\Gamma$-convergence of noncoercive functionals in $W^{1,p}
(\Omega )  $ (in chinese). {\sl Science in China Series A}, 1994, 24(3):
223-240

\item {9} Jian H.  Homogenization problems of parabolic minima. {\sl
Acta Math Appl Sinica}, 1996,  12: 318-327

\item{10}  Alt H,  Luckhaus S.  Quasilinear elliptic-parabolic differential
equations. {\sl Math Z}, 1982,  183: 311-341

\item{11} Tsutsumi M.  On solutions of some doubly nonlinear degenerate
parabolic equations with absorption.  {\sl J Math  Anal  Appl},
1988,  132: 187-212

\item{12}   Barbu V. Nonlinear semigroups and differential equations in
Banach space.  Amsterdam: Netherlands Press, 1976. 10-67

\end